\documentclass[a4,numbers,sort&compress]{elsarticle}%
\usepackage{amsfonts}
\usepackage{color}
\usepackage{amssymb,amsmath,amsthm,graphics,color}
\usepackage{graphics}
\usepackage{graphicx}
\usepackage{epsfig}
\usepackage{fancybox}
\usepackage{caption}
\usepackage{tikz,mathpazo}
\usepackage{multirow,booktabs}
\usepackage{amsmath}
\usepackage{amssymb}
\usepackage{geometry}
\usepackage[pagewise]{lineno}%
\providecommand{\U}[1]{\protect\rule{.1in}{.1in}}
\usetikzlibrary{shapes.geometric, arrows}
\usetikzlibrary{calc}
\journal{}

\newtheorem{theorem}{Theorem}[section]
\newtheorem{corollary}{Corollary}[section]
\newtheorem{remark}{Remark}[section]

\newtheorem{lemma}{Lemma}[section]
\begin{document}
\begin{frontmatter}
\title{Global stability of a time-delayed malaria model with  standard incidence rate}
\author[mymainaddress,mysecondaryaddress]{Songbai Guo}
\ead{guosongbai@bucea.edu.cn}
\author[mymainaddress]{Min He}
\author[mymainaddress]{Jing-An Cui\corref{mycorrespondingauthor}}
\ead{cuijingan@bucea.edu.cn}
\cortext[mycorrespondingauthor]{Corresponding author.}
\address[mymainaddress]{School of Science, Beijing University of Civil Engineering and Architecture, Beijing 102616, P.R. China}
\address[mysecondaryaddress]{Academy of Mathematics and Systems Science, Chinese Academy of Sciences,
Beijing 100190, P.R. China}
\begin{abstract}
A four-dimensional delay differential equations (DDEs) model of malaria with standard incidence rate is proposed. By utilizing the limiting system of the model and Lyapunov direct method, the global stability of equilibria of the model is obtained with respect to the basic reproduction number ${R}_{0}$. Specifically, it shows that the disease-free equilibrium ${E}^{0}$ is globally asymptotically stable (GAS) for ${R}_{0}<1$, and globally attractive (GA) for ${R}_{0}=1$, while the endemic equilibrium $E^{\ast}$ is GAS and ${E}^{0}$ is unstable for ${R}_{0}>1$. Especially, to obtain the global stability of the equilibrium $E^{\ast}$ for $R_{0}>1$, the weak persistence of the model is proved by some analysis techniques.
\end{abstract}
\end{frontmatter}


\section{Introduction}

Malaria is a disease caused by parasites and transmitted by the bites of
infected female anopheles mosquitoes, which can be life-threatening
\cite{World-1}. The World Health Organization \cite{World} estimated that there were 241 million malaria cases around the world in 2020, of which 627,000 deaths. Africa is the most affected region
with a serious influence on local economic development \cite{Snow}. Malaria is one of the most common infectious diseases in the world, bringing great pressure to the global control of infectious diseases \cite{Sachs}. The initial symptoms of malaria
(fever, headache and chills) generally arise 10 to 15 days after the bite of an infected mosquito \cite{World-1}, which indicates that malaria has a certain incubation period.

Since the emergence of malaria, disease prevention and prediction have been of
great importance. In this regard, mathematical modeling is a very important
method. Many scholars have established human-mosquito malaria models to
predict the development trend of the disease. In 1911, Ross \cite{Ross16} proposed a two-dimensional ordinary differential equations (ODEs) model of malaria and introduced the concept of threshold. In 1957, MacDonald \cite{Macdonald57} refined Ross's model and first defined the concept of basic reproduction number. The refined model is called the Ross-Macdonald model. In recent years, many scholars further revised the classical Ross-Macdonald model (see, e.g., \cite{Anderson1991,Rodriguez1997,Rodriguez2001}). Subsequent studies have extended the revised models to higher dimensions (see, e.g., \cite{Fatmawati21,Jin20,Mideo08,Ngwa06,Safan16,Wan09}). For example, in 2015, Safan and Ghazi \cite{Safan16} proposed a four-dimensional ODEs model of malaria with standard incidence rates, and analyzed the local and global stability of equilibria of the model. In malaria transmission, time delay can be considered due to the incubation period of the parasites, which means that individuals are not infectious until some time after they are infected \cite{Koutou19}. In 2008, Ruan et al. \cite{Ruan2008} first proposed the Ross-Macdonald model with time delay and studied the effect of time delay on the basic reproduction number. Afterwards, lots of scholars had established some DDEs models of malaria (see, e.g., \cite{Saker10,Li18,Ding19,Koutou19,Wu2021}). They mainly studied the local and global dynamics properties of these models. Consider that the standard incidence rate can be better to describe the transmission mechanism of malaria than the bilinear incidence rate (also see \cite{Li2020}). Thus, we propose a four-dimensional time-delayed malaria
model with standard incidence rate. Next, we will focus on the local and global stability of equilibria of the model with respect to the basic reproduction number $R_{0}$.

The remainder of this paper is organized as follows.
In Section 2, we first give model formulation, and
then obtain the well-posedness and dissipativeness of the
system. In Section 3, we explore the existence conditions of
equilibria of the system, and prove the local stability of equilibria with
respect to $R_{0}$. In
Section 4, in order to get the global stability of the endemic equilibrium for $R_{0}>1$, we verify the weak persistence of the system with a series of analysis techniques. In Section 5, we show that the global stability of the equilibria in terms of $R_{0}$ by using the limiting system of the model and
Lyapunov direct method.

\section{The model}

We divide humans and mosquitoes into four compartments, such as ${S}_{h}$:
susceptible humans, ${I}_{h}$: infected humans, ${S}_{v}$: susceptible
mosquitoes, ${I}_{v}$: infected mosquitoes. We use the positive parameters $\beta_{h}$, $\beta_{v}$, $\mu
_{h}$ and $\mu_{v}$ to represent the birth rates of humans and mosquitoes, and the death rates of humans and mosquitoes, respectively. The positive parameters $C_{vh}$ and $C_{hv}$ denote the infection rates of the infected mosquitoes biting susceptible humans and the susceptible mosquitoes biting the infected humans, respectively.
The nonnegative time delay $\tau$ is the incubation period of malaria. Then, the following model of malaria transmission is proposed:
\begin{equation}
\left\{
\begin{split}
\dot{S}_{h}(t)  &  ={\beta}_{h}-{C}_{vh}\frac{{I}_{v}(t)}{{N}_{v}(t)}{S}%
_{h}(t)-{\mu}_{h}{S}_{h}(t),\\
\dot{I}_{h}(t)  &  ={C}_{vh}\frac{{I}_{v}(t-\tau)}{{N}_{v}(t-\tau)}{S}%
_{h}(t-\tau)-{\mu}_{h}{I}_{h}(t),\\
\dot{S}_{v}(t)  &  ={\beta}_{v}-{C}_{hv}{I}_{h}(t){S}_{v}(t)-{\mu}_{v}{S}%
_{v}(t),\\
\dot{I}_{v}(t)  &  ={C}_{hv}{I}_{h}(t){S}_{v}(t)-{\mu}_{v}{I}_{v}(t),
\end{split}
\right.  \label{mod1}%
\end{equation}
where ${N}_{v}(t)={S}_{v}(t)+{I}_{v}(t)$.

The initial function of system \eqref{mod1} is
\[
\varphi=(\varphi_{1},\,\varphi_{2},\,\varphi_{3},\,\varphi_{4})^{T}\in{{C}_{+}%
}=\left\{  \varphi\in {C}={C}([-\tau,0],\,\mathbb{R}_{+}^{4}):\varphi_{3}%
(\theta)+\,\varphi_{4}(\theta)>0,\forall\theta\in\lbrack-\tau,0]\right\}
\]
with $\mathbb{R}_{+}=[0,\infty)$, where ${C}$ represents the Banach space composed of continuous mapping from $[-\tau,0]$ to $\mathbb{R}_{+}^{4}$ with the sup-norm. Next, we will discuss the well-posedness and the dissipativeness of system
\eqref{mod1} in ${{C}_{+}}$.

\begin{theorem}
\label{thm1} The solution $x(t)=(S_{h}(t),\,I_{h}(t),\,S_{v}(t),\,I_{v}%
(t))^{T}$ of system \eqref{mod1} through any $\varphi\in{{C}_{+}}$ exists, which is unique, nonnegative, and ultimately bounded on $[0,\infty)$.
\end{theorem}

\proof By using the similar proof of \cite[Proposition 2.1]{Guo18}, we first
can obtain the unique existence of the solution $x(t)$ and $x(t)\geq\mathbf{0}$ on
$[0,\infty)$. Hence, for $t\geq\tau$, it follows from system \eqref{mod1} that
\[
\dot{S}_{h}(t-\tau)+\dot{I}_{h}(t)={\beta}_{h}-{\mu}_{h}({S}_{h}(t-\tau
)+{I}_{h}(t)),\text{ }\dot{S}_{v}(t)+\dot{I}_{v}(t)=\beta_{v}-\mu_{v}%
(S_{v}(t)+I_{v}(t)),
\]
which yields
\[
\lim_{t\rightarrow\infty}(S_{h}(t-\tau)+I_{h}(t))=\frac{\beta_{h}}{\mu_{h}%
},\text{ }\lim_{t\rightarrow\infty}(S_{v}(t)+I_{v}(t))=\frac{\beta_{v}}%
{\mu_{v}}.
\]

\begin{remark}
It is not difficult to find that the solution $x(t)$ of system \eqref{mod1}
with any $\varphi\in{{C}_{+}}$ satisfies that $(S_{h}(t),S_{v}(t))^{T}\gg\mathbf{0}$ for $t>0$, and then ${{C}_{+}}$ is positively invariant for system \eqref{mod1}.
\end{remark}

\section{Local stability}

The disease-free equilibrium $E^{0}=(S_{h}^{0},0,S_{v}^{0},0)^{T}=(\frac{\beta_{h}}{\mu_{h}},0,\frac
{\beta_{v}}{\mu_{v}},0)^{T}$ of system \eqref{mod1} can be easily obtained. In order to get the existence of an endemic equilibrium (positive equilibrium) $E^{\ast}=(S_{h}^{\ast},I_{h}^{\ast},S_{v}^{\ast},I_{v}^{\ast
})^{T}$, the basic reproduction number $R_{0}=\sqrt{C_{vh}C_{hv}%
\beta_{h}/\mu_{h}\mu_{v}\mu_{h}}$ is first given by using the similar
method in \cite{Driessche02}.
\begin{lemma}
\label{thm2} System \eqref{mod1} has a unique endemic equilibrium $E^{\ast}$
if and only if $R_{0}>1$.
\end{lemma}

\proof By system \eqref{mod1}, the endemic equilibrium equations can be given
by
\begin{equation}
\left\{
\begin{split}
&  {\beta}_{h}-{C}_{vh}\frac{{I}_{v}^{\ast}}{S_{v}^{\ast}+I_{v}^{\ast}}{S}%
_{h}^{\ast}-{\mu}_{h}{S}_{h}^{\ast}=0,\\
&  {C}_{vh}\frac{{I}_{v}^{\ast}}{S_{v}^{\ast}+I_{v}^{\ast}}{S}_{h}^{\ast}%
-{\mu}_{h}{I}_{h}^{\ast}=0,\\
&  {\beta}_{v}-{C}_{hv}{I}_{h}^{\ast}{S}_{v}^{\ast}-{\mu}_{v}{S}_{v}^{\ast
}=0,\\
&  {C}_{hv}{I}_{h}^{\ast}{S}_{v}^{\ast}-{\mu}_{v}{I}_{v}^{\ast}=0.
\end{split}
\right.  \label{eq1}%
\end{equation}
Simplifying and sorting \eqref{eq1}, we have
\begin{equation}\label{eq00}
\left\{
\begin{split}
&  S_{h}^{\ast}=\frac{\beta_{h}}{C_{vh}\mu_{v}I_{v}^{\ast}/\beta_{v}+\mu_{h}%
},\\
&  I_{h}^{\ast}=\frac{\beta_{h}C_{vh}\mu_{v}I_{v}^{\ast}}{\beta_{v}\mu
_{h}\left(  C_{vh}\mu_{v}I_{v}^{\ast}/\beta_{v}+\mu_{h}\right)  },\\
&  S_{v}^{\ast}=\frac{\beta_{v}}{C_{hv}I_{h}^{\ast}+\mu_{v}},\\
&  I_{v}^{\ast}=\frac{\beta_{v}C_{hv}I_{h}^{\ast}}{\mu_{v}\left(  C_{hv}%
I_{h}^{\ast}+\mu_{v}\right)  }.
\end{split}
\right.
\end{equation}
It follows from the second and the fourth equations in \eqref{eq00} and $I_{h}^{*}>0$ that
\[
I_{h}^{*} \beta_{v} \mu_{h}\left[\frac{C_{v h} \mu_{v}}{\beta_{v}} \cdot \frac{C_{h v} I_{h}^{*} \beta_{v}}{\mu_{v}\left(C_{h v} I_{h}{ }^{*}+\mu_{v}\right)}+\mu_{h}\right]=C_{v h} \mu_{v} \beta_{h} \cdot \frac{C_{h v} I_{h}^{*} \beta_{v}}{\mu_{v}\left(C_{h v} I_{h}{ }^{*}+\mu_{v}\right)} .
\]
So, We can solve $I_{h}^{*}$ as
\[
I_{h}^{\ast}=\frac{C_{vh}C_{hv}\beta_{h}-\mu_{v}\mu_{h}\mu_{h}}{\mu_{h}%
C_{hv}\left(  C_{vh}+\mu_{h}\right)  }=\frac{\beta_{h}\mu_{v}\left(  R_{0}%
^{2}-1\right)  }{\beta_{h}C_{hv}+\mu_{v}\mu_{h}R_{0}^{2}}.
\]
Consequently, there hold
\[
{S}_{h}^{\ast}=\frac{{\beta}_{h}({C}_{hv}\frac{{\beta}_{h}}{\mu_{h}}+{\mu}%
_{v})}{{C}_{hv}\beta_{h}+{\mu}_{v}\mu_{h}R_{0}^{2}},\text{ }{S}_{v}^{\ast
}=\frac{{\beta}_{v}({C}_{vh}+{\mu}_{h})}{{C}_{vh}\mu_{v}+{\mu}_{v}\mu_{h}%
R_{0}^{2}},\text{ }{I}_{v}^{\ast}=\frac{{\beta}_{v}{\mu}_{h}(R_{0}^{2}-1)}%
{{C}_{vh}\mu_{v}+{\mu}_{v}\mu_{h}R_{0}^{2}}.
\]
Thus, the endemic equilibrium $E^{\ast}$ is uniquely obtained if and only
if $R_{0}>1.$

For the local asymptotic stability of the disease-free equilibrium $E^{0}$,
we have the following result.
\begin{theorem}
\label{thm3} For any $\tau\geq0$, $E^{0}$ is locally asymptotically stable (LAS)
if $R_{0}<1$ and unstable if $R_{0}>1$.
\end{theorem}

\proof By the direct calculation, the transcendental characteristic equation
of the linearized system of system \eqref{mod1} at $E^{0}$ is taken by%
\begin{equation}
(\lambda+\mu_{h})(\lambda+\mu_{v})[(\lambda+\mu_{h})(\lambda+\mu_{v}%
)-u_{1}u_{2}e^{-\lambda\tau}]=0, \label{equtc}%
\end{equation}
where $u_{1}={C}_{hv}\beta_{v}/\mu_{v},\,u_{2}=C_{vh}\beta_{h}\mu_{v}%
/\beta_{v}\mu_{h}.$ Obviously, equation \eqref{equtc} has two negative roots
$\lambda_{1}=-\mu_{v}$ and $\lambda_{2}=-\mu_{h}.$ Now, we consider the
following transcendental equation
\begin{equation}
G(\lambda)=\lambda^{2}+q_{1}\lambda+q_{2}+q_{3}e^{-\lambda\tau}=0,
\label{eq2203}%
\end{equation}
where $q_{1}=\mu_{h}+\mu_{v},\,q_{2}=\mu_{v}\mu_{h},\,q_{3}=-u_{1}u_{2}.$ Note
that
\[
q_{1}>0,\text{ }q_{2}+q_{3}=\mu_{v}\mu_{h}-{C}_{hv}\frac{\beta_{v}}{\mu_{v}%
}C_{vh}\frac{\beta_{h}\mu_{v}}{\beta_{v}\mu_{h}}=\mu_{v}\mu_{h}(1-R_{0}%
^{2})>0
\]
for $R_{0}<1.$ Thus, it is known from the Routh-Hurwitz criterion that each
root of equation \eqref{eq2203} has a negative real part for $R_{0}<1$ and
$\tau=0$.

For $R_{0}<1$ and $\tau>0$, assume that $G(\lambda)$ has the pure imaginary root
$\lambda=wi$ ($w\geq0$). Then it follows from equation \eqref{eq2203} that
\[
\left\{
\begin{split}
q_{2}-w^{2}  &  =-q_{3}\operatorname{cos}w\tau,\\
q_{1}w  &  =q_{3}\operatorname{sin}w\tau.
\end{split}
\right.
\]
Further, we have
\begin{equation}
w^{4}+\left(  q_{1}^{2}-2q_{2}\right)  w^{2}+q_{2}^{2}-q_{3}^{2}=0.
\label{equt1}%
\end{equation}
Clearly,
\[
q_{1}^{2}-2q_{2}=\mu_{v}^{2}+\mu_{h}^{2}>0,\text{ }q_{2}^{2}-q_{3}^{2}%
=(q_{2}-q_{3})(q_{2}+q_{3})>0
\]
for $R_{0}<1$ and $\tau>0.$ Thus, equation \eqref{equt1} is false, which means
that equation \eqref{eq2203} has no roots on the imaginary axis. Consequently,
each root of equation \eqref{eq2203} has a negative real part for $R_{0}<1$
and $\tau>0$. Therefore, $E^{0}$ is LAS for $R_{0}<1$ and $\tau\geq0$.

For $R_{0}>1$ and $\tau\geq0,$ it holds that
\[
G(0)=q_{2}+q_{3}=\mu_{v}\mu_{h}(1-R_{0}^{2})<0,\text{ }\lim\limits_{\lambda
\rightarrow\infty}G(\lambda)=\infty.\text{ }%
\]
In consequence, equation \eqref{eq2203} must admit a positive root, which
implies that $E^{0}$ is unstable for $R_{0}>1$ and $\tau\geq0.$

For the local asymptotic stability of the endemic equilibrium $E^{\ast}$,
we have the following conclusion.
\begin{theorem}
\label{thm4} If $R_{0}>1$, then $E^{\ast}$ is LAS for any $\tau\geq0$.
\end{theorem}

\proof By the direct calculation, the transcendental characteristic equation of the linearized system of system \eqref{mod1} at $E^{\ast}$ is given by
\begin{equation}
\left\vert
\begin{array}
[c]{cccc}%
\lambda+{C}_{vh}\frac{{I}_{v}^{\ast}}{{N}_{v}^{\ast}}+\mu_{h} & 0 & -{C}%
_{vh}\frac{{I}_{v}^{\ast}{S}_{h}^{\ast}}{{{N}_{v}^{\ast}}^{2}} & {C}_{vh}%
\frac{{S}_{v}^{\ast}{S}_{h}^{\ast}}{{{N}_{v}^{\ast}}^{2}}\\
-{C}_{vh}\frac{{I}_{v}^{\ast}}{{N}_{v}^{\ast}}e^{-\lambda\tau} & \lambda
+\mu_{h} & {C}_{vh}\frac{{I}_{v}^{\ast}{S}_{h}^{\ast}}{{{N}_{v}^{\ast}}^{2}%
}e^{-\lambda\tau} & -{C}_{vh}\frac{{S}_{v}^{\ast}{S}_{h}^{\ast}}{{{N}%
_{v}^{\ast}}^{2}}e^{-\lambda\tau}\\
0 & {C}_{hv}{S}_{v}^{\ast} & \lambda+{C}_{hv}{I}_{h}^{\ast}+\mu_{v} & 0\\
0 & -{C}_{hv}{S}_{v}^{\ast} & -{C}_{hv}{I}_{h}^{\ast} & \lambda+\mu_{v}%
\end{array}
\right\vert =0. \label{eq3}%
\end{equation}
Let $m_{1}={C}_{vh}{I}_{v}^{\ast}/{N}_{v}^{\ast},\,m_{2}={C}_{vh}%
{I}_{v}^{\ast}{S}_{h}^{\ast}/{{N}_{v}^{\ast}}^{2},\,m_{3}={C}%
_{vh}{S}_{v}^{\ast}{S}_{h}^{\ast}/{{N}_{v}^{\ast}}^{2},\,m_{4}%
={C}_{hv}{S}_{v}^{\ast},\,m_{5}={C}_{hv}{I}_{h}^{\ast}$. Then equation
\eqref{eq3} becomes%

\begin{equation}
(\lambda+\mu_{v})(\lambda+\mu_{h})[(\lambda+\mu_{h}+m_{1})(\lambda+\mu
_{v}+m_{5})-m_{4}(m_{3}+m_{2})e^{-\lambda\tau}]=0. \label{eq4}%
\end{equation}
Obviously, equation \eqref{eq4} has characteristic roots $\lambda_{1}=-\mu
_{v}<0,\,\lambda_{2}=-\mu_{h}<0.$ Now, we consider the following
transcendental equation
\begin{equation}
G(\lambda)=\lambda^{2}+p_{1}\lambda+p_{2}+p_{3}e^{-\lambda\tau}, \label{eq5}%
\end{equation}
where $p_{1}=\mu_{h}+m_{1}+\mu_{v}+m_{5},\,p_{2}=(\mu_{h}+m_{1})(\mu_{v}%
+m_{5}),\,p_{3}=-m_{4}(m_{3}+m_{2}).$

For $R_{0}>1$ and $\tau=0$, we have%
\[%
\begin{split}
p_{1}  &  =\mu_{h}+m_{1}+\mu_{v}+m_{5}>0,\\
p_{2}+p_{3}  &  =\left(  \mu_{h}+m_{1}\right)  \left(  \mu_{v}+m_{5}\right)
-m_{4}\left(  m_{3}+m_{2}\right) \\
&  =C_{hv}\left(  \mu_{h}+C_{vh}\right)  I_{h}^{\ast}+C_{vh}\frac{\mu_{v}%
}{\beta_{v}}\left(  \mu_{v}+C_{hv}\frac{\beta_{h}}{\mu_{h}}\right)
I_{v}^{\ast}+\mu_{h}\mu_{v}\left(  1-R_{0}^{2}\right) \\
&  =\left(  R_{0}+1\right)  \mu_{v}\mu_{h}\left(  R_{0}-1\right)  >0,
\end{split}
\]
where ${I}_{h}^{\ast}={\beta}_{h}{\mu}_{v}(R_{0}^{2}-1)/({C}_{hv}\beta
_{h}+{\mu}_{v}\mu_{h}R_{0}^{2}),\,{I}_{v}^{\ast}={\beta}_{v}{\mu}_{h}%
(R_{0}^{2}-1)/({C}_{vh}\mu_{v}+{\mu}_{v}\mu_{h}R_{0}^{2})$ and $R_{0}%
^{2}=C_{vh}C_{hv}\beta_{h}/\mu_{h}\mu_{v}\mu_{h}$ are used. It follows from
the Routh-Hurwitz criterion that each root of equation \eqref{eq5} has a
negative real part.

For $R_{0}>1$ and $\tau>0$, if $G(\lambda)$ has the root $\lambda=wi$
($w\geq0$), then by equation \eqref{eq5}, we have
\begin{equation*}
\left\{
\begin{split}
p_{2}-w^{2}  &  =-p_{3}\operatorname{cos}w\tau\\
p_{1}w  &  =p_{3}\operatorname{sin}w\tau.
\end{split}
\right.
\end{equation*}
Consequently, it holds
\begin{equation}
w^{4}+\left(  p_{1}^{2}-2p_{2}\right)  w^{2}+p_{2}^{2}-p_{3}^{2}=0.
\label{equt2}%
\end{equation}
Note that%
\begin{align*}
p_{1}^{2}-2p_{2}  &  =\left(  \mu_{h}+m_{1}\right)  ^{2}+\left(  \mu_{v}%
+m_{5}\right)  ^{2}+2m_{4}\left(  m_{3}+m_{2}\right)  >0,\\
p_{2}^{2}-p_{3}^{2}  &  =\left(  p_{2}-p_{3}\right)  \left(  p_{2}%
+p_{3}\right)  >0.
\end{align*}
Accordingly, equation \eqref{equt2} is false, which hints that equation
\eqref{eq5} has no roots on the imaginary axis. Hence, each root of equation
\eqref{eq5} has a negative real part for $R_{0}>1$ and $\tau>0$. Thus,
$E^{\ast}$ is LAS for $R_{0}>1$ and $\tau\geq0$.

\section{Weak persistence}

To get the global stability of $E^{\ast}$, we need to study the weak persistence of
system \eqref{mod1}. Let $D=\{\varphi\in{C}_{+}:\varphi_{2}(0)>0\}$ and
\[
x(t)=(S_{h}(t),I_{h}(t),S_{v}(t),I_{v}(t))^{T}%
\]
be the solution of system \eqref{mod1} through any $\varphi\in D.$ It is not difficult to find that $D$ is a positive invariant set with respect to system \eqref{mod1} and $x(t)\gg\mathbf{0}$ for $t>0$. Consequently, we investigate the weak persistence of system \eqref{mod1} in $D$. Following the definition of weak persistence in \cite{Butler86}, system \eqref{mod1} is called weakly persistent if $\limsup_{t\rightarrow\infty}\mathcal{u}(t)>0$, $\mathcal{u}=S_{h},I_{h},S_{v}%
,I_{v}$ for any $\varphi\in D$. Define $x_{t}=\left(  S_{ht},I_{ht},S_{vt},I_{vt}\right)^{T}\in {C}_{+}$ for $t\geq0$ as $x_{t}(\theta)=x(t+\theta),$ $\theta\in\lbrack-\tau,0].$ Then $x_{t}$ is the solution of system \eqref{mod1} through $\varphi$.
Motivated by the persistence approach in
\cite{Guo22}, we have the following results.

\begin{lemma}
\label{lem1} If ${R}_{0}>1$, $\theta\in(0,1)$ , and $\sideset{}{}\limsup
_{t\rightarrow\infty}I_{h}(t)\leq\theta I_{h}^{*}$, then there hold
\[
\liminf_{t\rightarrow\infty}S_{v}(t)\ge\bar{S}_{v}\equiv\frac{\beta_{v}}{
\theta C_{hv} I_{h}^{*}+\mu_{v}}>{S}_{v}^{*},~~\liminf_{t\rightarrow\infty
}S_{h}(t)\ge\bar{S}_{h}\equiv\frac{\beta_{h}}{C_{vh}(1-\frac{\bar{S}_{v}%
}{S_{v}^{0}})+\mu_{h}}>{S}_{h}^{*} .
\]

\end{lemma}

\proof According to the endemic equilibrium equations, we have $\bar{S}_{v}>{S}_{v}^{\ast},\,\bar{S}_{h}>{S}_{h}^{\ast}.$ For any $\varepsilon>1$, there exists ${\rho}={\rho}(\varphi,\varepsilon)>0$ such that for $t>0$, it follows
\[
I_{h}(t)<\varepsilon\theta I_{h}^{\ast}.
\]
Therefore, for $t>{\rho}$, we get
\[
\dot{S}_{v}(t)=\beta_{v}-C_{vh}I_{h}(t)S_{v}(t)-\mu_{v}S_{v}(t)>\beta
_{v}-S_{v}(t)\left(  \varepsilon\theta C_{hv}I_{h}^{\ast}+\mu_{v}\right)  .
\]
As a result, there holds
\[
\liminf_{t\rightarrow\infty}S_{v}(t)\geq\frac{\beta_{v}}{\varepsilon\theta
C_{hv}I_{h}^{\ast}+\mu_{v}}.
\]
Let $\varepsilon\rightarrow1^{+}$, we obtain
\[
\liminf_{t\rightarrow\infty}S_{v}(t)\geq\bar{S}_{v}.
\]
It follows from system \eqref{mod1} that $\dot{N_{v}}(t)=\beta_{v}-\mu
_{v}N_{v}(t)$, which yields
\begin{equation}
\label{eq14}\lim_{t\rightarrow\infty}N_{v}(t)=S_{v}^{0}.
\end{equation}
Hence, we have
\[
\limsup_{t\rightarrow\infty}\frac{I_{v}(t)}{N_{v}(t)}=\limsup_{t\rightarrow
\infty}\left(  1-\frac{S_{v}(t)}{N_{v}(t)}\right)  =1-\liminf_{t\rightarrow
\infty}\frac{S_{v}(t)}{N_{v}(t)}\leq1-\frac{\bar{S}_{v}}{S_{v}^{0}}.
\]
For any $\epsilon\in(0,1)$, there is $\tilde{\rho}=\tilde{\rho}(\varphi
,\epsilon)>0$, such that for $t>\tilde{\rho}$,
\[
\frac{I_{v}(t)}{{N}_{v}(t)}<1-\frac{\epsilon\bar{S}_{v}}{{S}_{v}^{0}}.
\]
Hence, for $t>\tilde{\rho}$, we have
\[
\dot{S}_{h}(t)=\beta_{h}-C_{vh}\frac{I_{v}(t)}{N(t)}S_{h}(t)-\mu_{h}%
S_{h}(t)>\beta_{h}-S_{h}(t)\left[  C_{vh}\left(  1-\frac{\epsilon\bar{S}_{v}%
}{S_{v}^{0}}\right)  +\mu_{h}\right]  .
\]
Consequently, there holds
\[
\liminf_{t\rightarrow\infty}S_{h}(t)\geq\frac{\beta_{h}}{C_{vh}(1-\frac
{\epsilon\bar{S}_{v}}{S_{v}^{0}})+\mu_{h}}.
\]
Let $\epsilon\rightarrow1^{-}$, we obtain
\[
\liminf_{t\rightarrow\infty}S_{h}(t)\geq\bar{S}_{h}.
\]

\begin{theorem}
\label{thm7} If $R_{0}>1$ and $\theta\in(0,1)$, then $\limsup_{t\rightarrow
\infty}I_{h}(t)>\theta I_{h}^{\ast}$.
\end{theorem}

\proof We prove the conclusion by contradiction. Assume $\limsup_{t\rightarrow\infty}I_{h}%
(t)\leq\theta I_{h}^{\ast}$. Then according to Lemma \ref{lem1}), for any $\varepsilon\in(0,\varepsilon_{0})$, there exists an $\varepsilon_{0}>0$ such that
\begin{equation*}
\frac{\bar{S}_{h}}{S_{v}^{0}+\varepsilon}>\frac{{S}_{h}^{\ast}}{S_{v}^{0}}.
\end{equation*}
 By Lemma \ref{lem1} and
\eqref{eq14}, we get that for any $\varepsilon\in(0,\varepsilon_{0}),$ there is a $\mathcal{T}\equiv \mathcal{T}\left(  \varepsilon,\varphi\right)>0$ such that
\[
\frac{S_{h}(t)}{N_{v}(t)}>\frac{\bar{S}_{h}}{S_{v}^{0}+\varepsilon}%
,\,S_{v}(t)>S_{v}^{\ast}
\]
for all $t\geq \mathcal{T}$.
Next, we define the following functional
\[
L(\varphi)=\varphi_{2}+\frac{C_{vh}\bar{S}_{h}}{(S_{v}^{0}+\varepsilon)\mu
_{v}}\varphi_{4}+{C}_{vh}\int_{-\tau}^{0}\frac{\varphi_{4}(\theta)}%
{\varphi_{3}(\theta)+\varphi_{4}(\theta)}\varphi_{1}(\theta)d\theta.
\]
Clearly, $L$ is continus on $\mathbb{R}_{+}^{4}$ and then ${L}(x_{t})$ is bounded. Now, we calculate the derivative of $L$ along the solution $x_{t}$ for $t\geq \mathcal{T}$ as follows
\begin{align*}
\dot{L}(x_{t})  &  \geq\left[  \frac{C_{vh}C_{hv}\bar{S}_{h}}{(S_{v}%
^{0}+\varepsilon)\mu_{v}}S_{v}(t)-\mu_{h}\right]  I_{h}(t)\\
&  >\mu_{h}\left[  \frac{{S}_{v}^{\ast}\bar{S}_{h}}{(S_{v}^{0}+\varepsilon
)S_{h}^{0}}R_{0}^{2}-1\right]  I_{h}(t).
\end{align*}

Set
\[
\bar{I}_{h}=\min_{\theta\in[-\tau,0]}I_{h}(\mathcal{T}+\tau+\theta
),\,\Lambda=\min\left\{  \bar{I}_{h},\frac{\mu_{v}I_{v}(\mathcal{T})}{C_{hv}{S}%
_{v}^{\ast}}\right\}  \text{.}%
\]
Now, we start by verifying $I_{h}(t)\geq\Lambda$ for $t\geq \mathcal{T}.$ Otherwise, there is a $\mathcal{T}_{0}\geq0$ such that $I_{h}(t)\geq\Lambda$ for $t\in\lbrack \mathcal{T},\mathcal{T}+\mathcal{T}_{0}+\tau]$, $I_{h}(\mathcal{T}+\mathcal{T}_{0}+\tau)=\Lambda$ and $\dot{I}_{h}(\mathcal{T}+\mathcal{T}_{0}%
+\tau)\leq0$. For $t\in\lbrack \mathcal{T},\mathcal{T}+\mathcal{T}_{0}+\tau]$, we can obtain
\begin{equation}
\dot{I}_{v}(t)={C}_{hv}{I}_{h}(t){S}_{v}(t)-{\mu}_{v}{I}_{v}(t)\geq{C}%
_{hv}\Lambda{S}_{v}^{\ast}-{\mu}_{v}{I}_{v}(t). \label{eq16}%
\end{equation}
By \eqref{eq16}, for $t\in\lbrack \mathcal{T},\mathcal{T}+\mathcal{T}_{0}+\tau]$, we have
\[
I_{v}(t)\geq\dfrac{{C}_{hv}\Lambda{S}_{v}^{\ast}}{{\mu}_{v}}+\left(
I_{v}(T)-\dfrac{{C}_{hv}\Lambda{S}_{v}^{\ast}}{{\mu}_{v}}\right)
e^{uT-ut}\geq\dfrac{{C}_{hv}\Lambda{S}_{v}^{\ast}}{{\mu}_{v}}%
\]
for $t\in\lbrack \mathcal{T},\mathcal{T}+\mathcal{T}_{0}+\tau]$. Note that $R_{0}^{2}=S_{v}^{0}S_{h}^{0}%
/S_{v}^{\ast}S_{h}^{\ast}$, we have
\[
\frac{{S}_{v}^{\ast}\bar{S}_{h}}{(S_{v}^{0}+\varepsilon)S_{h}^{0}}R_{0}%
^{2}-1>\frac{{S}_{v}^{\ast}{S}_{h}^{\ast}}{S_{v}^{0}S_{h}^{0}}R_{0}^{2}-1=0.
\]
As a consequence, we conclude that
\begin{align*}
\dot{I}_{h}(\mathcal{T}+\mathcal{T}_{0}+\tau)  &  ={C}_{vh}\frac{{I}_{v}(\mathcal{T}+\mathcal{T}_{0})}{{N}_{v}%
(T+T_{0})}{S}_{h}(\mathcal{T}+\mathcal{T}_{0})-{\mu}_{h}{I}_{h}(\mathcal{T}+\mathcal{T}_{0}+\tau)\\
&  >\Lambda\mu_{h}\left[  \frac{{S}_{v}^{\ast}\bar{S}_{h}}{(S_{v}%
^{0}+\varepsilon)S_{h}^{0}}R_{0}^{2}-1\right]>0,
\end{align*}
which contradicts $\dot{I}_{h}(\mathcal{T}+\mathcal{T}_{0}+\tau)\leq0$. Accordingly, ${I}_{h}(t)\geq\Lambda$ for $t\geq \mathcal{T}\text{.}$ Therefore, for $t\geq T$, we have
\[
\dot{L}(x_{t})\geq\Lambda\mu_{h}\left[  \frac{{S}_{v}^{\ast}\bar{S}_{h}%
}{(S_{v}^{0}+\varepsilon)S_{h}^{0}}R_{0}^{2}-1\right]  >0,
\]
which yields $L(x_{t})\rightarrow\infty$ as $t\rightarrow\infty\text{.}$ As a consequence,
this is a contradiction to the boundedness of $L(x_{t})$.

From Theorem \ref{thm7}, it is not difficult to obtain the following corollary.

\begin{corollary}
\label{coro1} If $R_{0}>1$, then system \eqref{mod1} is weakly persistent for $\tau\geq0$.
\end{corollary}

\section{Global stability}

In this section, we will prove the global stability of $E^{0}\,$and
$E^{\ast}$ with respect to $R_{0}$. To do this, it follows from \eqref{eq14}
that system \eqref{mod1} has the limiting system as follows:%
\begin{equation}
\left\{
\begin{split}
\dot{S}_{h}(t)  &  ={\beta}_{h}-\frac{{C}_{v}h{I}_{v}(t){S}_{h}(t)}{S_{v}^{0}%
}-{\mu}_{h}{S}_{h}(t),\\
\dot{I}_{h}(t)  &  =\frac{{C}_{vh}{I}_{v}(t-\tau){S}_{h}(t-\tau)}{S_{v}^{0}%
}-{\mu}_{h}{I}_{h}(t),\\
\dot{S}_{v}(t)  &  ={\beta}_{v}-{C}_{hv}{I}_{h}(t){S}_{v}(t)-{\mu}_{v}{S}%
_{v}(t),\\
\dot{I}_{v}(t)  &  ={C}_{hv}{I}_{h}(t){S}_{v}(t)-{\mu}_{v}{I}_{v}(t).
\end{split}
\right.  \label{mod2}%
\end{equation}
Following the proof of Theorem \ref{thm1}, the solution $y(t)=\left(
S_{h}(t),I_{h}(t),S_{v}(t),I_{v}(t)\right)  ^{T}$ of system \eqref{mod2} with
any $\psi\in {C}_{+}$ exists, and it is unique, nonnegative, and ultimately bounded on $[0,\infty)$. Obviously, $E^{0}$ and $E^{\ast}$ are the equilibria of system \eqref{mod2}. In addition, ${C}_{+}$ is positively invariant for system \eqref{mod2} and $(S_{h}(t),S_{v}(t))^{T}\gg\mathbf{0}$ for $t>0.$ Let us define $y_{t}=\left(  S_{ht},I_{ht},S_{vt},I_{vt}\right)  ^{T}\in {C}_{+}$ for $t\geq0$ as $y_{t}(\theta)=y(t+\theta),$ $\theta\in\lbrack-\tau,0].$ Then we can obtain the following the global stability result of the disease-free equilibrium
$E^{0}$ of system \eqref{mod1}.

\begin{theorem}
\label{thm5} For any $\tau\geq0$, $E^{0}$ is GAS
for $R_{0}<1$ and GA for $R_{0}=1$ in ${C}_{+}$.
\end{theorem}

\proof Firstly, it follows from Theorem \ref{thm3} that $E^{0}$ is LAS for $R_{0}<1$. We thus just require to show that $E^{0}$ is GA for $R_{0}\leq1$. Let $x_{t}$ be the solution of system \eqref{mod1} through any $\varphi\in {C}_{+}$ and $y_{t}$ be the solution of system \eqref{mod2} with any $\psi\in {C}_{+}$. By Theorem \ref{thm1}, we have that $x_{t}$ is bounded. Hence, we can obtain that $\omega(\varphi)\subseteq {C}_{+}$ is compact, where $\omega(\varphi)$ is the $\omega$-limit set of $\varphi$ for system \eqref{mod1}. To verify that $E^{0}$ is GA, we just require to prove $\omega(\varphi)=\{E^{0}\}.$

Next, we define a functional $V$ on
$
\Omega_{1}=\left\{  \psi\in {C}:\psi_{1}(0)>0,\psi_{3}(0)>0\right\}  \subseteq
{C}_{+}%
$
as follows%
\begin{equation}
V\left(  \psi\right)  =V_{1}(\psi(0))+\int_{-\tau}^{0}\frac{\mu_{v}}{\beta
_{v}}C_{vh}\psi_{4}(\theta)\psi_{1}(\theta)d\theta, \label{eq20}%
\end{equation}
where
\[
V_{1}(\psi(0))=S_{h}^{0}\left(  \frac{\psi_{1}(0)}{S_{h}^{0}}%
-1-\ln\frac{\psi_{1}(0)}{S_{h}^{0}}\right)  +\psi_{2}(0)+\frac{\mu_{v}\mu_{h}%
}{C_{hv}\beta_{v}}S_{v}^{0}\left(  \frac{\psi_{3}(0)}{S_{v}^{0}}-1-\ln
\frac{\psi_{3}(0)}{S_{v}^{0}}\right)  +\frac{\mu_{v}\mu_{h}}{C_{hv}\beta_{v}%
}\psi_{4}(0).
\]
Clearly, $V_{1}(\psi(0))\leq V\left(  \psi\right)  $ and $V_{1}$ is positive
definite with respect to $E^{0}$ on $\Omega_{1}.$ We calculate the derivative of $V$ along the solution $y_{t}$ for $t\geq 1$ as follows
\begin{align}
\dot{V}\left(  y_{t}\right)  =  &  \left(  1-\frac{S_{h}^{0}}{S_{h}}\right)
\dot{S}_{h}+\dot{I}_{h}+\frac{\mu_{v}\mu_{h}}{C_{hv}\beta_{v}}\left(
1-\frac{S_{v}^{0}}{S_{v}}\right)  \dot{S}_{v}+\frac{\mu_{v}\mu_{h}}%
{C_{hv}\beta_{v}}\dot{I}_{v}+\frac{\mu_{v}}{\beta_{v}}C_{vh}I_{v}%
S_{h}\nonumber\\
&  -\frac{\mu_{v}}{\beta_{v}}C_{vh}I_{v}(t-\tau)S_{h}(t-\tau)\nonumber\\
=  &  -\frac{\mu_{h}\left(  S_{h}^{0}-S_{h}\right)  ^{2}}{S_{h}}-\frac{\mu
_{v}\mu_{h}}{C_{hv}\beta_{v}}\frac{\mu_{v}\left(  S_{v}^{0}-S_{v}\right)
^{2}}{S_{v}}+\left(  \frac{\beta_{h}\mu_{v}}{\mu_{h}\beta_{v}}C_{vh}-\frac
{\mu_{v}\mu_{h}}{C_{hv}\beta_{v}}\mu_{v}\right)  I_{v}\nonumber\\
=  &  -\frac{\mu_{h}\left(  S_{h}^{0}-S_{h}\right)  ^{2}}{S_{h}}-\frac{\mu
_{v}\mu_{h}}{C_{hv}\beta_{v}}\frac{\mu_{v}\left(  S_{v}^{0}-S_{v}\right)
^{2}}{S_{v}}+\frac{\mu_{v}\mu_{v}\mu_{h}\left(  R_{0}^{2}-1\right)  }%
{C_{hv}\beta_{v}}I_{v}\leq0. \label{eq9}%
\end{align}
By \eqref{eq20} and \eqref{eq9}, it follows from \cite[Corollary 3.3]{Guo}
that $E^{0}$ is uniformly stable for system \eqref{mod2}. Furthermore, we can
have that $S_{h}(t)$ and $S_{v}(t)$ are persistent, that is, there is some
$\epsilon=\epsilon(\psi)>0$ such that $\liminf\nolimits_{t\rightarrow\infty
}S_{i}(t)>\epsilon,$ $i=h,v$. Hence, we obtain $\omega(\psi)\subseteq
\Omega_{1},$ where $\omega(\psi)$ is the $\omega$-limit set of $\psi$ in regard to system \eqref{mod2}. Consequetly, $V$ is a Lyapunov functional on $\{y_{t}:t\geq1\}\subseteq\Omega_{1}$. \cite[Corollary 2.1]{Guo} implies that $\,\dot{V}=0$ on $\omega(\psi).$

Let $y_{t}$ be the solution of system \eqref{mod2} with any $\phi\in
\omega(\psi)$. By the invariance of $\omega(\psi)$, we can obtain that $y_{t}$
$\in\omega(\psi)$ for $t\in\mathbb{R}$. From \eqref{eq9}, it follows that
$S_{h}(t)=S_{h}^{0}$ and $S_{v}(t)=S_{v}^{0}$ for $t\in\mathbb{R}$. The first
and the third equations of system \eqref{mod2} hint that $I_{h}%
(t)=I_{v}(t)=0$ for $t\in\mathbb{R}$. Thus, we have $\omega(\psi)=\{E^{0}\}$,
which means that $W^{s}\left(  E^{0}\right)  ={C}_{+}$, where $W^{s}\left(
E^{0}\right)$ is the stable set of $E_{0}$ for system \eqref{mod2}.
Obviously, $\omega(\varphi)\cap W^{s}\left(  E^{0}\right)  \neq\varnothing$.
In consequence, it follows from \cite[Theorem 4.1]{Thieme92} that
$\omega(\varphi)=\{E_{0}\}$.

For global stability of the endemic equilibrium $E^{\ast}$ of system \eqref{mod1}, we have
the following result.
\begin{theorem}
\label{thm6} If $R_{0}>1$, then $E^{\ast}$ is GAS for any $\tau\geq0$ in $D$.
\end{theorem}

\proof By Theorem \ref{thm4}, we have that $E^{\ast}$ is LAS for $R_{0}>1$.
Hence, we only need to show that $E^{\ast}$ is GA for $R_{0}>1$. Let $x_{t}$
and $y_{t}$ be the solution of system \eqref{mod1} with any $\varphi\in D$ and
any $\psi\in D,$ respectively. Note that $D$ is a positive invariant set for
system \eqref{mod2}, and $x(t)\gg\mathbf{0}$ for $t>0$. Let ${\Omega
}_{2}=\{\psi\in{C}:\psi(0)\gg\mathbf{0}\}$. Then $\Omega_{2}{\subseteq D}.$ To prove that $E^{\ast}$ is GA for $R_{0}>1$, we just require to verify $\omega(\varphi)=\{E^{\ast}\}.$

Next, we define a functional $V$ on ${\Omega}_{2}$ as follows
\begin{equation}
V\left(  \psi\right)  =V_{2}{}(\psi(0))+\mu_{h}I_{h}^{\ast}\int_{-\tau}%
^{0}\left(  \frac{\mu_{v}C_{vh}\psi_{4}(\theta)\psi_{1}(\theta)}{\beta_{v}%
\mu_{h}I_{h}^{\ast}}-1-\ln\frac{\mu_{v}C_{vh}\psi_{4}(\theta)\psi_{1}(\theta
)}{\beta_{v}\mu_{h}I_{h}^{\ast}}\right)  d\theta, \label{eq10}%
\end{equation}
where%
\begin{align}
V_{2}{}  &  (\psi(0))=\left(  \psi_{1}(0)-S_{h}^{\ast}-S_{h}^{\ast}\ln
\frac{\psi_{1}(0)}{S_{h}^{\ast}}\right)  +\left(  \psi_{2}(0)-I_{h}^{\ast
}-I_{h}^{\ast}\ln\frac{\psi_{2}(0)}{I_{h}^{\ast}}\right) \nonumber\\
&  +\frac{\mu_{h}I_{h}^{\ast}}{\mu_{v}I_{v}^{\ast}}\left(  \psi_{3}%
(0)-S_{v}^{\ast}-S_{v}^{\ast}\ln\frac{\psi_{3}(0)}{S_{v}^{\ast}}\right)
+\frac{\mu_{h}I_{h}^{\ast}}{\mu_{v}I_{v}^{\ast}}\left(  \psi_{4}%
(0)-I_{v}^{\ast}-I_{v}^{\ast}\ln\frac{\psi_{4}(0)}{I_{v}^{\ast}}\right)
.\nonumber
\end{align}
It easily follows that $V_{2}(\psi(0))\leq V\left(  \psi\right)  $ and $V_{2}$
is positive definite with respect to $E^{\ast}$ on $\Omega_{2}.$ We calculate the derivative of $V$ along the solution $y_{t}$ for $t\geq 1$ as follows
\begin{align*}
\dot{V}\left(  y_{t}\right)  ={}  &  \left(  1-\frac{S_{h}^{\ast}}{S_{h}%
}\right)  \dot{S}_{h}+\left(  1-\frac{I_{h}^{\ast}}{I_{h}}\right)  \dot{I}%
_{h}+\frac{\mu_{h}I_{h}^{\ast}}{\mu_{v}I_{v}^{\ast}}\left(  1-\frac
{S_{v}^{\ast}}{S_{v}}\right)  \dot{S}_{v}+\frac{\mu_{h}I_{h}^{\ast}}{\mu
_{v}I_{v}^{\ast}}\left(  1-\frac{I_{v}^{\ast}}{I_{v}}\right)  \dot{I}%
_{v}\nonumber\\
&  +\frac{\mu_{v}}{\beta_{v}}C_{vh}I_{v}S_{h}-\frac{\mu_{v}}{\beta_{v}}%
C_{vh}I_{v}(t-\tau)S_{h}(t-\tau)+\mu_{h}I_{h}^{\ast}\ln\frac{I_{v}%
(t-\tau)S_{h}(t-\tau)}{I_{v}S_{h}}\nonumber\\
=  &  -\frac{\mu_{h}\left(  S_{h}-S_{h}^{\ast}\right)  ^{2}}{S_{h}}-\frac
{\mu_{h}I_{h}^{\ast}}{\mu_{v}I_{v}^{\ast}}\frac{\mu_{v}\left(  S_{v}%
-S_{v}^{\ast}\right)  ^{2}}{S_{v}}+4\mu_{h}I_{h}{}^{\ast}-\mu_{h}I_{h}{}%
^{\ast}\frac{S_{h}{}^{\ast}}{S_{h}}-\mu_{h}I_{h}{}^{\ast}\frac{I_{h}^{\ast}%
}{I_{h}}\frac{I_{v}(t-\tau)S_{h}(t-\tau)}{I_{v}{}^{\ast}S_{h}{}^{\ast}%
}\nonumber\\
&  -\mu_{h}I_{h}{}^{\ast}\frac{S_{v}^{\ast}}{S_{v}}-\mu_{h}I_{h}{}^{\ast}%
\frac{I_{v}{}^{\ast}}{I_{v}}\frac{I_{h}S_{v}}{I_{h}{}^{\ast}S_{v}{}^{\ast}%
}+\mu_{h}I_{h}{}^{\ast}\ln\frac{I_{v}(t-\tau)S_{h}(t-\tau)}{I_{v}S_{h}}.
\end{align*}
Using the following equality
\[
\ln\frac{I_{v}(t-\tau)S_{h}(t-\tau)}{I_{v}S_{h}}=\ln\frac{S_{h}^{\ast}}{S_{h}%
}+\ln\frac{I_{h}^{\ast}I_{v}(t-\tau)S_{h}(t-\tau)}{I_{h}I_{v}^{\ast}S_{h}^{\ast}}+\ln\frac{S_{v}^{\ast}}{S_{v}}
+\ln\frac{I_{v}^{\ast}I_{h}S_{v}}{I_{v}I_{h}^{\ast}S_{v}^{\ast}},
\]
we can obtain that
\begin{align}
\dot{V}\left(  y_{t}\right)  ={}  &  -\frac{\mu_{h}\left(  S_{h}-S_{h}^{\ast
}\right)  ^{2}}{S_{h}}-\frac{\mu_{h}I_{h}^{\ast}}{\mu_{v}I_{v}^{\ast}}%
\frac{\mu_{v}\left(  S_{v}-S_{v}^{\ast}\right)  ^{2}}{S_{v}}+\mu_{h}%
I_{h}^{\ast}F\left(  \frac{S_{h}^{\ast}}{S_{h}}\right) \nonumber\\
&  \quad+\mu_{h}I_{h}^{\ast}F\left(  \frac{I_{h}^{\ast}}{I_{h}}\frac
{I_{v}(t-\tau)S_{h}(t-\tau)}{I_{v}^{\ast}S_{h}^{\ast}}\right)  +\mu_{h}%
I_{h}^{\ast}F\left(  \frac{S_{v}^{\ast}}{S_{v}}\right)  +\mu_{h}I_{h}^{\ast
}F\left(  \frac{I_{v}^{\ast}}{I_{v}}\frac{I_{h}S_{v}}{I_{h}^{\ast}S_{v}^{\ast
}}\right)  , \label{eq12}%
\end{align}
where $F(x)=1-x+\ln x$, $x>0$. Since $F(x)$ is a nonnegative difinite function
for $x>0,$ $\dot{V}\left(  y_{t}\right) \leq0$ for $t\geq1.$ Thus, by \eqref{eq10},
\eqref{eq12} and \cite[Corollary 3.3]{Guo}, we have that $E^{\ast}$ is
uniformly stable for system \eqref{mod2}. Moreover, we can obtain $\omega(\psi)\subseteq\Omega_{2}$, which yields that $V$ is a Lyapunov functional on $\{y_{t}:t\geq1\}\subseteq\Omega_{2}$. It follows from
\cite[Corollary 2.1]{Guo} that $\,\dot{V}=0$ on $\omega(\psi).$

Let $y_{t}$ is the solution of system \eqref{mod2} with any $\phi\in
\omega(\psi)$. The invariance of $\omega(\psi)$ hints that $y_{t}$ $\in
\omega(\psi)$ for $t\in\mathbb{R}$. It follows from \eqref{eq12} that
\[
S_{h}(t)=S_{h}^{\ast},\text{ }S_{v}(t)=S_{v}^{\ast},\text{ }I_{h}^{\ast}%
I_{v}(t)=I_{v}^{\ast}I_{h}(t)
\]
for $t\in\mathbb{R}$. Thus, by the fourth equation of system \eqref{mod2}, we
have
\[
I_{h}^{\ast}\dot{I}_{v}=I_{h}^{\ast}C_{hv}I_{h}S_{v}-I_{h}^{\ast}\mu_{v}%
I_{v}=I_{h}^{\ast}C_{hv}I_{h}S_{v}^{\ast}-\mu_{v}I_{h}I_{v}^{\ast}%
=(C_{hv}I_{h}^{\ast}S_{v}^{\ast}-\mu_{v}I_{v}^{\ast})I_{h}=0,
\]
where $C_{hv}I_{h}^{\ast}S_{v}^{\ast}=\mu_{v}I_{v}^{\ast}$ is used.
Consequently, $I_{v}(t)$ is a constant function and then $I_{h}(t)$ is also a
constant function. In consequence, $y(t)\gg\mathbf{0}$ is an endemic
equilibrium of system \eqref{mod2} for $R_{0}>1$. Considering that the
endemic equilibrium is unique, we thus have $y(t)=E^{\ast}$. Hence, it holds
that $\omega(\psi)=E^{\ast}$, which means that $W^{s}\left(  E^{\ast}\right)
=D$, where $W^{s}\left( E^{\ast}\right)  $ is the stable set of $E^{\ast}$ for
system \eqref{mod2}. From Theorem \ref{thm7}, we have $\omega
(\varphi)\cap W^{s}\left(  E^{\ast}\right)  \neq\varnothing$.
Thus, \cite[Theorem 4.1]{Thieme92} indicates that $\omega(\varphi)=\{E^{\ast}\}$.

\section*{Acknowledgements}

This work is supported in part by the National Natural Science Foundation of
China (Nos. 11901027 and 11871093), and the China Postdoctoral Science
Foundation (No. 2021M703426), the Pyramid Talent Training Project of BUCEA
(No. JDYC20200327).

\section*{References}


\end{document}